\documentclass[11pt]{article}
\usepackage{amssymb}
\newcommand{\be}{\begin{equation}}
\newcommand{\ee}{\end{equation}}
\newcommand{\ra}{\rightarrow}

\newtheorem{theorem}{Theorem}
\newtheorem{lemma}{Lemma}
\begin{document}
\begin{center}
{\bf Efficiency requires innovation}
\end{center}
\begin{center}
Abram Kagan \\
Department  of Mathematics, University of Maryland\\
 College Park, MD 20742, USA
\end{center}
\begin{center}
{\bf Abstract}
\end{center}
\noindent
In estimation a parameter $\theta\in{\mathbb R}$ from a sample $(x_1,\ldots,x_n)$ from a population $P_{\theta}$ a simple way of incorporating a new observation $x_{n+1}$ into an estimator $\tilde\theta_{n}=\tilde\theta_{n}(x_1,\ldots,x_n)$ is transforming $\tilde\theta_n$ to what we call the
{\it jackknife extension} $\tilde\theta_{n+1}^{(e)}=\tilde\theta_{n+1}^{(e)}(x_1,\ldots,x_n,x_{n+1})$,
\[\tilde\theta_{n+1}^{(e)}=\{\tilde\theta_n (x_1 ,\ldots,x_n)+
\tilde\theta_n (x_{n+1},x_2 ,\ldots,x_n)+\ldots+\tilde\theta_n (x_1 ,\ldots,x_{n-1},x_{n+1})\}/(n+1).\]
Though $\tilde\theta_{n+1}^{(e)}$ lacks an innovation the statistician could expect from a larger data set, it is still better than $\tilde\theta_n$,
\[{\rm var}(\tilde\theta_{n+1}^{(e)})\leq\frac{n}{n+1} {\rm var}(\tilde\theta_n).\]
However, an estimator obtained by jackknife extension for all $n$ is asymptotically efficient only for samples from exponential families.
For a general $P_{\theta}$, asymptotically efficient estimators require innovation when a new observation is added to the data.\\
Some examples illustrate the concept.
\newpage
\section{Introduction}
Let $\tilde\theta_n=\tilde\theta_{n}(x_1,\ldots,x_n)$ be an estimator of $\theta$ based on sample of size $n$ from a population $P_\theta$ with $\theta\in{\mathbb R}$ as a parameter. If another observation $x_{n+1}$ is added to the data, a simple way of incorporating it in the
existing estimator is by what we call the {\it jackknife extension},
\be
\tilde\theta_{n+1}^{(e)}=\tilde\theta_{n+1}^{(e)}(x_1,\ldots,x_n,x_{n+1})=(\tilde\theta_{n,1}+\ldots+\tilde\theta_{n,n+1})/(n+1)
\ee
where
\[\tilde\theta_{n,i} =\tilde\theta_n (x_1,\ldots,x_{i-1},x_{i+1},\ldots, x_{n}),\:i=2,\ldots,n,\:\tilde\theta_{n,n+1} =\tilde\theta_n (x_1,\ldots,x_n, x_{n}).\]
Plainly, $E_{\theta}(\tilde\theta_{n+1}^{(e)}) =E_{\theta}(\tilde\theta_{n})$ and if $\tilde\theta_n$ is symmetric in its arguments (as is usually the case) the jackknife extension is symmetric in
$x_1,\ldots,x_n,x_{n+1}$.\\
\\
 If ${\rm var}_{\theta}(\tilde\theta_{n})<\infty$, then not only ${\rm var}_{\theta}(\tilde\theta_{n+1}^{(e)})<{\rm var}_{\theta}(\tilde\theta_{n})$ but a stronger inequality holds:
\be
(n+1){\rm var}_{\theta}(\tilde\theta_{n+1}^{(e)})\leq n{\rm var}_{\theta}(\tilde\theta_n).
\ee
The inequality (2) is a direct corollary of a special case of the so called {\it variance drop lemma} due to (Artstein {\it et al.},2004).
\begin{lemma} Let $X_1,\ldots,X_n,X_{n+1}$ be independent identically distributed random variables and $\psi(X_1,\ldots,X_n)$ a function with
$E(|\psi(X_1,\ldots,X_{n+1})|^2)<\infty$. Set
\[\psi_1 =\psi(X_2,\ldots,X_{n+1}),\psi_i
=\psi(X_1,\ldots,X_{i-1}, X_{i+1},\ldots,X_{n+1}),\:i=2,\ldots,n+1.\]
Then
\be
{\rm var}(\sum_{1}^{n+1}\psi_i)\leq n\sum_{1}^{n+1}{\rm var}(\psi_i).
\ee
\end{lemma}
Note that with $n+1$ instead of $n$ on the right hand side of (5), the inequality becomes a trivial corollary of
\[(\sum_{1}^{n+1} a_i)^2\leq (n+1)\sum_{1}^{n+1} a_{i}^2\]
holding for any numbers $a_1,\ldots,a_{n+1}$.\\
 For an extension of the variance drop lemma see (Madiman an Barron, 2007).\\
Suppose that starting with $n=m$ and $\tilde\theta_{m}(x_1,\ldots,x_m)$, the statistician constructs the jackknife extension
$\tilde\theta_{m+1}^{(e)}$ of $\tilde\theta_{m}(x_1,\ldots,x_m)$, then the jackknife extension $\tilde\theta_{m+2}^{(e)}$ of
$\tilde\theta_{m+1}^{(e)}$ and so on. One can easily see that for $n\geq 2m$ the estimator $\tilde\theta_{n}^{(e)}(x_1,\ldots,x_n)$ thus obtained is a classical $U$-statistic with the kernel $\tilde\theta_{m}(x_1,\ldots,x_m)$:
\be
\tilde\theta_{n}^{(e)}(x_1,\ldots,x_n)=\frac{1}{{n\choose m}}\sum_{1\leq i_{1}\leq\ldots\leq i_{m}\leq n} \tilde\theta_{m}(x_{i_1},\ldots,x_{i_m}).
\ee
Hoeffding initiated studying $U$-statistics back in 1948.
The variance of $\tilde\theta_{n}^{(e)}$ can be explicitly expressed in terms of $\tilde\theta_{m}$. Set
\[\tilde\theta_{m|k}(x_1,\ldots,x_k)=E_{\theta}\{\tilde\theta_{m}(X_1,\ldots,X_m) |X_1 =x_1,\ldots,X_k =x_k\}.\]
The following formula due to Hoeffding (1948) expresses ${\rm var}(\tilde\theta_{n}^{(e)})$ via
$v_k (\theta)={\rm var}(\tilde\theta_{m|k}(X_1,\ldots,X_k)),\:k=1,\ldots,m$:
\be
{\rm var}_{\theta}(\tilde\theta_{n}^{(e)})=\frac{1}{{n\choose m}}\sum_{k=1}^{m}{m\choose k}{n-m \choose m-k}v_k (\theta).
\ee
\section{Main result}
Assume that the distributions $P_{\theta}$ are given by differentiable in $\theta$ density (with respect to
a measure $\mu$)  $p(x;\theta)$ with the Fisher information
\[I(\theta)=\int(\frac{\partial\log p(x;\theta)}{\partial\theta})^2p(x;\theta)d\mu(x)\]
well defined and finite.\\
If $E_{\theta}(\tilde\theta_m)=\gamma(\theta)$, by Cram\'{e}r-Rao inequality
\[{\rm var}_{\theta}(\tilde\theta_{m})\geq\frac{|\gamma'(\theta)|^2}{mI(\theta)}.\]
In particular, if $\tilde\theta_m$ is an unbiased estimator of $\theta$, ${\rm var}_{\theta}(\tilde\theta_m)\geq\frac{1}{mI(\theta)}$.\\
Furthermore, if ${\rm var}_{\theta}(\tilde\theta_m)<\infty$, then ${\rm var}_{\theta}(\tilde\theta_n^{(e)})<\infty$ for all $n>m$ and the following lemma holds.
\begin{lemma}(Hoeffding 1948). As $n\rightarrow\infty$, $\sqrt{n}(\tilde\theta_{n}^{(e)} -\gamma(\theta))$
is asymptotically normal $N(0,m^2 v_{1}(\theta))$.
\end{lemma}
Due to Cram\'{e}r-Rao inequality, for any unbiased estimator $\tilde\gamma_n$ of $\gamma(\theta)$ based on a sample $(x_1,\ldots,x_n)$ from a population with Fisher information $I(\theta)$,
\be
{\rm var}_{\theta}(\gamma_n)\geq\frac{|\gamma'(\theta)|^2}{nI(\theta)}.
\ee
Combining (6) with Lemma 2 leads to a formula for the asymptotic efficiency of $\tilde\theta_n^{(e)}$:
\be
{\rm aseff}(\tilde\theta_{n}^{(e)})=\frac{|\gamma'(\theta)|^2 /I(\theta)}{m^2 v_{1}(\theta)}.
\ee
\begin{lemma} Let $X$ be a random element, $X\sim p(x;\theta)$ with finite Fisher information $I(\theta)$. If $h(X)$ is a (scalar valued) function with $E_{\theta}(h(X))=\mu(\theta)$ differebtiable and ${\rm var}_{\theta}(h(X)=\sigma^{2}(\theta)<\infty$, then
\be
I(\theta)\geq\frac{|\mu'(\theta)|^2}{\sigma^{2}(\theta)}.
\ee
\end{lemma}
{\it Proof}. Take the projection of the Fisher score
$J(X;\theta)=(p'(x;\theta)/p(x;\theta)$ into the subspace span${1,h(X)}$ of the Hilbert space of functions with $g(X)$ with $E_{\theta}(|g(X)|^2)<\infty$:
\be
\hat E_{\theta}\{(J(X;\theta)|1, h(X)\}=\hat J(X;\theta)=a(\theta)(h(X)-\mu(\theta)).
\ee
Multiplying both sides by $h(X)-\mu(\theta)$ and taking the expectations results in $a(\theta)=\mu'(\theta)/\sigma^{2}(\theta)$ due to the property
$E_{\theta}(J(X;\theta)h(X))=\mu'(\theta)$ of the Fisher score. Hence
\[I(\theta)={\rm var}_{\theta}(J(X;\theta))\geq{\rm var}_{\theta}(\hat J(X;\theta))=\frac{|\mu'(\theta)|^2}{\sigma^{2}(\theta)}\]
which is exactly (8). The equality sign in (8) is attained if and only if with
$P_{\theta}$-probability one the relation
\be
\frac{p'(x;\theta)}{p(x;\theta)}=a(\theta)(h(x)-\gamma(\theta))
\ee
holds for $a(\theta)$. \\
From $E_{\theta}(\tilde\theta_{m|1})=\gamma(\theta), v_{1}={\rm var}(\tilde\theta_{m|1})$ and (8)
one gets
\be
{\rm aseff}(\tilde\theta_{n}^{(e)})\leq 1/m^2.
\ee
Thus, a necessary condition for the asymptotic efficiency of $\tilde\theta_{n}^{(e)}$ is $m=1$ and by virtue of (4)
\be
\tilde\theta_{n}^{(e)} (x_1,\ldots,x_n)=(h(x_1) +\ldots +h(x_n))/n
\ee
for some $h(x)$ with $E_{\theta}\{h(X)\}=\gamma(\theta)$.\\
From Lemma 3 the estimator (12) is an asymptotically efficient estimator of $\gamma(\theta)$ if and only if
the relation (10) holds implying that the family is exponential,
\be
p(x;\theta)=\exp \{A(\theta)h(x)+B(\theta)+g(x)\}
\ee
where the functions in the exponent are such that $E_{\theta}(h(X))=\gamma(\theta)$.\\
From (8) one can see that that the maximum likelihood equation for $\theta$ based on a sample $(x_1,\ldots,x_n)$ from population (13) is
\be
(h(x_1) +\ldots +h(x_n))/n = \gamma(\theta)
\ee
and
\[\tilde\theta_{n}^{(e)} (x_1,\ldots,x_n)=(h(x_1) +\ldots +h(x_n))/n\]

as the maximum likelihood estimator of $\gamma(\theta)$ is asymptotically efficient.\\
We summarize the above as a theorem.
\begin{theorem} Under the regularity type conditions of the theory of maximum likelihood estimators, the jackknife extension estimators are asymptotically efficient if and only if they are arithmetic means based on samples from exponential families.
\end{theorem}
\section{Some examples}
The jackknife extension lacks innovation. A jackknife extension estimator based on $(x_1,\ldots,x_{n+1})$ differs from the estimator based on $(x_,\ldots,x_n)$ only by the sample sample size. In a sense, it is an extensive vs. intensive use of the data when the main factor is quantity vs. quality.\\
\\
Nonparametric estimators of population characteristics such as the empirical distribution function, the sample mean and variance are jackknife extensions. Their main goal is  to be universal rather than optimal for individual populations. An interesting statistic is the sample median $\tilde\mu_n =\tilde\mu_n (x_1,\ldots,x_n)$ constructed from a sample from a continuous population. Without loss in generality, one may assume
\[x_1 <\ldots<x_n.\]
For $n=2m+1,\:\tilde\mu_n =x_{m+1}$. If $x_{n+1}<x_1$ or $x_{n+1}>x_n$, one can easily see that
\be
\tilde\mu_{n+1}^{(e)}=(x'_{m+1} + x'_{m+2})/2
\ee
where $x'_{m+1}$ and $x'_{m+2}$ are the $(m+1)$st and $(m+2)$nd elements of the sample $(x_1,\ldots,x_{n+1})$.
Thus, the median of a sample of an even size is a jackknife extension though one should keep in mind that the definitions of the
median in samples of even and odd size are different and it is not clear if the inequality (2) holds.\\
\\
For $n=2m$ the jackknife extension of $\tilde\mu_{n}=(x_m + x_{m+1})/2$ is not $x'_{m+1}$. Let us start with simple cases of
$m=2$ and $m=3$. In the first case,
\[\tilde\mu_{5}^{(e)}=\frac{1.5x'_2 +2x'_3 +1.5x'_4}{5}\]
is a weighted average of $x'_3$ and its nearest neighbors. The same holds in the second case, with $x'_4$ instead of $x'_3$ and different
weights:
\[\tilde\mu_{7}^{(e)}=\frac{2x'_3 +3x'_4 +2x'_5}{7}.\]
It seems likely that the extrapolation to an arbitrary $n=2m$ will result in
\be
\tilde\mu_{n+1}^{(e)}=\frac{((m+1)/2) x'_{m} +m x'_{m+1} +((m+1)/2) x'_{m+2}}{n+1}.
\ee
Though (16) is a reasonable estimator of the median, it is not clear how it behaves for $n=2m$ in small and large samples compared to the standard $\tilde\mu_{n+1} =x'_{m+1}$.
\section{References}
Artstein, S., Ball, K. M., Barthe, F., Naor, A. (2004).
Solution of Shannon's problem on the monotonicity of entropy.
{\it J. Amer. Math. Soc.,} {\bf 17,}  975--982.\\
\\
Hoeffding, W. (1948).
A class of statistics with asymptotically normal distribution.
{\it Ann. Math. Stat.,} {\bf 19,} 293--325.\\
\\
Kagan, A. M., Yu, T., Barron, A., Madiman, M. (2011).
Contribution to the theory of Pitman estimators.
{\it J. Math. Sci.},{\bf 199,} 2, 202-214.\\
\\
 Madiman, M., Barron, A. (2007).
 Generalized Entropy Power Inequalities and Monotonicity Properties of Information.
 {\it IEEE Transactions on Information Theory,} {\bf 53,} 2317--2329.
\end{document}